\documentclass[11pt,reqno]{amsart}

\usepackage[T1]{fontenc}
\usepackage[utf8]{inputenc}
\usepackage{amsmath,amsfonts,amssymb,amsthm}
\usepackage[margin=1.1in]{geometry}
\usepackage[colorlinks=true,linkcolor=blue,citecolor=blue,urlcolor=blue]{hyperref}

\renewcommand{\le}{\leqslant}
\renewcommand{\ge}{\geqslant}

\newcommand{\PP}{\mathsf{P}}
\newcommand{\EE}{\mathsf{E}}
\newcommand{\DD}{\mathsf{D}}
\newcommand{\bC}[2]{C_{#1}^{#2}}
\newcommand{\qb}[2]{\genfrac{[}{]}{0pt}{}{#1}{#2}_{q}}
\newcommand{\I}{\mathcal I}
\DeclareMathOperator{\stdeg}{deg}

\theoremstyle{plain}
\newtheorem{theorem}{Theorem}
\newtheorem*{thmA}{Theorem A}
\newtheorem{corollary}{Corollary}
\newtheorem{lemma}{Lemma}
\newtheorem{proposition}{Proposition}
\theoremstyle{definition}
\newtheorem{remark}{Remark}

\numberwithin{equation}{section}

\begin{document}

\title[Threshold probabilities for stability of the independence number]
      {On Threshold Probabilities for the Stability of the\\
       Independence Number of a Random Subgraph}

\author{V.~A.~Pokhachevskiy}
\address{Lomonosov Moscow State University, Moscow, Russia}
\email{pokhachevskiy@gmail.com}

\author{A.~M.~Raigorodskii}
\address{Moscow Institute of Physics and Technology (National Research
University), Dolgoprudny, Moscow Region, Russia;
Lomonosov Moscow State University, Moscow, Russia;
Caucasus Mathematical Center, Adyghe State University, Maikop, Russia;
Buryat State University, Institute for Mathematics and Informatics,
Ulan-Ude, Russia}
\email{mraigor@yandex.ru}


\date{\today}

\begin{abstract}
For an arbitrary sequence of graphs $G_n$, we consider a random
spanning subgraph $G_{n,p}$ in which each edge of the graph $G_n$ survives
independently of the others with probability $p$. We obtain a general upper bound
for the independence number $\alpha(G_{n,p})$, from which, for an arbitrary
probability $p=p(n)$, we derive sufficient conditions for asymptotic
stability, i.e., for the equality $\alpha(G_{n,p})=(1+o(1))\alpha(G_n)$
asymptotically almost surely, as well as upper bounds for the independence
number of the random subgraph, including unconditional ones; for $p=1/2$ they recover
A.\,M.~Raigorodskii's results of 2017 and 2026.
Under additional conditions on the structure of the
independent sets of the graph $G_n$, we find, up to a constant factor, the
threshold probabilities for both asymptotic and exact stability (the equality
$\alpha(G_{n,p})=\alpha(G_n)$). Special cases include results on the graphs
$G(n,r,<s)$ for growing $r$ and $s$, including Kneser graphs for
growing $r$; we also give
examples of complete multipartite graphs, graphs on the blocks of Steiner systems, and
$q$-analogues of Kneser graphs.

\end{abstract}

\maketitle

\section{Introduction}
\label{sec:intro}

Recall that for a simple graph $G=(V,E)$ (undirected, without loops or multiple
edges) the independence number is the quantity $\alpha(G)$ equal to the
maximum number of pairwise nonadjacent vertices of the graph; the sets of pairwise
nonadjacent vertices themselves are called independent. The classical binomial model
of a random graph is the Erd\H{o}s--R\'enyi model $G(n,p)$, i.e., a random graph on
$n$ labeled vertices in which each of the $\bC n2$ possible edges
is drawn independently of the others with probability $p$ (see~\cite{Boll,JLR,Kolchin,RaiBook}).
The independence number $\alpha(G(n,p))$ is well studied; in particular, for many
$p=p(n)$, asymptotically almost surely (hereafter, a.a.s.) one has
$\alpha(G(n,p))\sim2\log_{1/(1-p)}(np)$ \cite{Boll,JLR}. This statement
involves two asymptotic relations. A sequence of events $\mathcal A_n$
occurs a.a.s.\ if $\PP(\mathcal A_n)\to1$ as $n\to\infty$; for
random variables $X_n$ and a positive function $f(n)$, the notation
``a.a.s.\ $X_n\sim f(n)$'' (or ``a.a.s.\ $X_n=(1+o(1))f(n)$'') means
that for every fixed $\varepsilon>0$
\begin{displaymath}
\PP\bigl((1-\varepsilon)f(n)\le X_n\le(1+\varepsilon)f(n)\bigr)\to1,
\qquad n\to\infty,
\end{displaymath}
i.e., $X_n/f(n)\to1$ in probability. One-sided
bounds of the form ``a.a.s.\ $X_n\le(1+o(1))f(n)$'' are understood similarly (for every $\varepsilon>0$,
a.a.s.\ $X_n\le(1+\varepsilon)f(n)$), and the notation ``a.a.s.\ $X_n\asymp f(n)$''
means that for some constants $0<c\le C$ a.a.s.\
$cf(n)\le X_n\le Cf(n)$.

The graph $G(n,p)$ is a random spanning subgraph of the complete graph $K_n$. Following
\cite{Rai17}, consider an arbitrary sequence of graphs
$G_n=(V_n,E_n)$, $n=1,2,\dots$, and a random element $G_{n,p}$ taking values in the
set of all spanning subgraphs of the graph $G_n$ and having the distribution
\begin{displaymath}
\PP\bigl(G_{n,p}=(V_n,E)\bigr)=p^{|E|}(1-p)^{|E_n|-|E|},\qquad E\subseteq E_n .
\end{displaymath}
Deleting edges does not decrease the independence number; therefore, one always has
$\alpha(G_{n,p})\ge\alpha(G_n)$. We say that for the independence numbers of the graphs
$G_n$ \emph{asymptotic stability} holds if a.a.s.\
$\alpha(G_{n,p})\sim\alpha(G_n)$, and \emph{exact stability} holds if a.a.s.\
$\alpha(G_{n,p})=\alpha(G_n)$ (in \cite{Rai17} it is called strong, and in
\cite{RaiSin} absolute stability). For example, there is no stability for the complete graph:
$\alpha(K_n)=1$, whereas $\alpha(G(n,p))\to\infty$ a.a.s.\ for every
constant $p<1$.

Random subgraphs of \emph{distance graphs} have been studied most extensively; the vertices
of such graphs are the $r$-element subsets of the set $[n]=\{1,\dots,n\}$ and whose edges
are determined by the cardinalities of the intersections. For Kneser graphs (the edges are pairs
of disjoint sets), asymptotic stability for $p=1/2$ and constant $r\ge2$ was proved in
\cite{BGPR1,BGPR2}; then B.~Bollob\'as, B.~Narayanan, and A.\,M.~Raigorodskii
\cite{BNR} found, for $2\le r=o(n^{1/3})$, the threshold for exact stability up
to a factor $1+o(1)$, and
these results were substantially extended in \cite{Kup16,DasTran,DevlinKahn,BKL}.
For the graphs $G(n,r,<s)$, in which edges join sets intersecting
in fewer than $s$ elements (for $s=1$ this is the Kneser graph), stability
was studied in \cite{Pyad15,OgRai,KOR,RaiKaras,PR}; in \cite{PR2} the threshold
for exact stability was found up to a constant factor for growing
$r=r(n)$ and $s=s(n)$. For the graphs $G(n,r,s)$, in which edges join
sets intersecting in exactly $s$ elements, M.\,M.~Pyaderkin
\cite{Pyad19} found, for constant $r$ and $s$ with $r>2s+1$, the threshold for asymptotic
stability. For a detailed survey of these and many other results, see
\cite{RaiSin}.

For arbitrary sequences of graphs, the first general result on
asymptotic stability was obtained in \cite{Rai17}. We state it in its original formulation.

\begin{thmA}[\cite{Rai17}]
Let $G_n=(V_n,E_n)$, $n\in\mathbb N$, be a sequence of graphs. Let
$N_n=|V_n|$, $\alpha_n=\alpha(G_n)$. Let $\gamma_n$ be the maximum
number of vertices of the graph $G_n$ that are nonadjacent to both vertices of a given edge.
Suppose that the quantities $N_n,\alpha_n,\gamma_n$ tend monotonically to
infinity and that there exists a function $\beta_n$ such that the following
asymptotic conditions hold:
1.~$\beta_n>\gamma_n$ and $\beta_n=o(\alpha_n)$;
2.~$\log_2N_n=o(\alpha_n/\beta_n)$;
3.~$\log_2N_n=o(\beta_n-\gamma_n)$.
Then for $p=\frac12$ asymptotic stability holds for the independence numbers of the graphs $G_n$.
\end{thmA}

Recently, in \cite{Rai26}, Theorem~A was generalized in another direction: for $p=1/2$,
instead of asymptotic stability, the independence number is allowed to grow by a factor of
$\lambda_n\ge1$; namely, under conditions of the form 1--3 with an additional
parameter $\lambda_n$, it is proved that a.a.s.\ $\alpha(G_{n,1/2})\le(1+o(1))\alpha_n\lambda_n$
(Theorem~3 of \cite{Rai26}). In the same paper, the unconditional upper bound
$\alpha(G_{n,1/2})\le(1+o(1))\cdot2\alpha_n\log_2N_n$ and lower bounds for the
independence number in terms of the maximum vertex degree and in terms of families of pairwise
nonadjacent subgraphs are obtained; the latter apply, in particular, to the graphs
$G(n,2s+1,s)$.

The key idea of the proof of Theorem~A goes back to M.\,M.~Pyaderkin's argument for
Kneser graphs \cite{BGPR2} (it is stated as a separate lemma in
\cite[Subsection~3.3.10, Lemma~2]{RaiSin}): if two adjacent vertices have few common
non-neighbors, then Tur\'an's bound for the number of edges in a vertex set can be
substantially strengthened. In the graphs $G(n,r,<s)$ with $n$ sufficiently large
compared with $r$ and $s$, the independent sets of maximum cardinality are the
\emph{stars}, i.e., the families of all $r$-sets containing a fixed
$s$-set; the survey
\cite[Subsection~3.3.10]{RaiSin} notes that for arbitrary graphs
``the notion of a \emph{star} loses its meaning'' and that the problem of finding more flexible conditions
for asymptotic stability remains open.

In the present paper, the role of stars is played by all independent sets of maximum
cardinality. Instead of bounds on the cardinality of nontrivial intersecting families,
which for Kneser graphs and the graphs $G(n,r,<s)$ are provided by the
Hilton--Milner \cite{HM} and Ahlswede--Khachatrian \cite{AK2} theorems, we use the
quantity $\psi(G_n)$, the maximum cardinality of an independent set not
contained in any independent set of maximum cardinality. We obtain the following results.
\begin{itemize}
\item We prove a general bound for the probability of the event $\alpha(G_{n,p})\ge k$
(Theorem~\ref{th:core}), which allows us to extend Theorem~A,
Theorem~3 of \cite{Rai26}, and the unconditional bound $2\alpha_n\log_2N_n$ to
the case of an arbitrary probability $p=p(n)$, moreover with an explicit bound on the deviation
$\alpha(G_{n,p})-\alpha(G_n)$ (Corollaries~\ref{sled:quant}--\ref{sled:uncond}).
\item We prove a Tur\'an-type bound (Theorem~\ref{th:turan}) and find the order of growth of the quantity $\alpha(G_{n,p})$ for small $p$
(Theorem~\ref{th:order} and Corollary~\ref{sled:order}); along the way, we obtain an analogue of
the lower bound from \cite{Rai26} in terms of pairwise nonadjacent subgraphs for
arbitrary $p$ (Remark~\ref{zam:lower}).
\item Under conditions on the quantity $\psi(G_n)$ and some other characteristics
of the graphs, we find, up to a constant factor, the threshold for asymptotic
stability (Theorem~\ref{th:asthr}) and the threshold for exact stability
(Theorem~\ref{th:exthr}).
\end{itemize}
Theorem~\ref{th:exthr} contains as a special case the main result
of \cite{PR2} on the graphs $G(n,r,<s)$; other examples are complete multipartite
graphs, the construction from \cite{Rai17}, graphs on the blocks of Steiner systems (they
include the Steiner network example from \cite{Rai17}), and $q$-analogues of Kneser
graphs.

The paper is organized as follows. In Section~\ref{sec:res} we introduce notation and
state the results, in Section~\ref{sec:ex} we analyze examples,
Section~\ref{sec:proofs} contains the proofs of the theorems, and Section~\ref{sec:concl}
contains concluding remarks.

\section{Notation and statements of the results}
\label{sec:res}

In what follows, $G_n=(V_n,E_n)$ is a sequence of graphs, $N=N_n=|V_n|$,
$\alpha=\alpha_n=\alpha(G_n)$, $p=p(n)\in[0,1]$; as a rule, we omit the index $n$
and write $G=G_n$, $G_p=G_{n,p}$. All asymptotic relations
refer to $n\to\infty$. We assume that the graph $G$ has at least one
edge, i.e., $N>\alpha$. For $D\subseteq V$, we denote by $e(D)$ the number of
edges of the graph $G$ with both endpoints in $D$. Throughout, $C_n^r$ denotes the binomial
coefficient $\binom{n}{r}$. We will need the following
characteristics of the graph $G$.

1.~The quantity $\gamma=\gamma(G)$ from Theorem~A: the maximum, over all edges
$\{x,y\}\in E$, of the number of vertices $z\in V\setminus\{x,y\}$ adjacent neither to $x$ nor
to $y$ (\emph{common non-neighbors} of the endpoints of the edge).

2.~The family $\I=\I(G)$ of all independent sets of cardinality $\alpha$ and the number
$m=|\I|$ of such sets.

3.~The quantity $\psi=\psi(G)$, the maximum cardinality of an independent set
not contained in any $I\in\I$ (if there are no such sets, then $\psi:=0$).
Since every independent set is contained in an inclusion-maximal
independent set, $\psi$ is the maximum cardinality of an inclusion-maximal
independent set that is not of maximum cardinality; in
particular, $\psi<\alpha$.

4.~The \emph{degree} of a set $D\subseteq V$ with respect to $\I$:
$\stdeg D:=\max_{I\in\I}|D\cap I|$ (for the graphs $G(n,r,<s)$, this is the degree of a family in the sense of \cite{OgRai}: the maximum, over
all stars, of the cardinality of the intersection of the family with a star).

5.~The quantity $\mu=\mu(G)$, the minimum, over all $I\in\I$ and all $v\in V\setminus
I$, of the number of vertices of the set $I$ adjacent to $v$, and $\delta:=1-\mu/\alpha$. Since
$I\cup\{v\}$ is not independent, $\mu\ge1$; clearly, $\mu\le\alpha$.

6.~The quantity $\theta=\theta(G)=\max\{|I\cap I'|\colon I,I'\in\I,\ I\ne I'\}$
(if $m=1$, then $\theta:=0$).

7.~The quantity $\Lambda=\Lambda(G)=\ln\bigl(m(N-\alpha)\bigr)$, the logarithm of the number of
pairs $(I,v)$, where $I\in\I$, $v\in V\setminus I$.

We note two simple facts. First, $m\le\bC N\alpha\le(eN/\alpha)^\alpha$,
i.e.,
\begin{equation}
\label{eq:lnm}
\ln m\le\alpha\ln\frac{eN}{\alpha}.
\end{equation}
Second, the following monotonicity property in $p$ holds.

\begin{remark}
\label{zam:monot}
For $p'\ge p$, the random graph $G_{p'}$ can be realized as the union of
$G_p$ and a random subgraph $G_{q}$ independent of it, where
$1-p'=(1-p)(1-q)$, and then every set independent in
$G_{p'}$ is independent in $G_p$. Therefore, the probabilities of the events
$\{\alpha(G_p)\ge k\}$ are nonincreasing in $p$, and the probability of the event
``$\alpha(G_p)=\alpha$ and $\I(G_p)=\I(G)$'' is nondecreasing in $p$ (the sets from
$\I(G)$ are independent in $G_p$ for every $p$). Consequently, statements of the form ``for $p\ge p_0$ (respectively, $p\le p_0$)
a.a.s.\ the event $\mathcal A$ occurs'', where $\mathcal A$ is one of the
above events with parameters independent of $p$, need to be proved only
for $p=p_0$.
\end{remark}

The proofs of all statements of this section are given in
Section~\ref{sec:proofs}; after each theorem and each corollary we indicate in which subsection
of Section~\ref{sec:proofs} it is proved.

\subsection{Arbitrary edge probability}

We begin with a general bound from which all the results of this subsection are derived.

\begin{theorem}
\label{th:core}
Let $p\in(0,1)$, let $k$ be an integer, $\alpha<k\le N$, and let $\beta$ be a
real number. If $\beta>\gamma$, then
\begin{equation}
\label{eq:core}
\PP\bigl(\alpha(G_p)\ge k\bigr)\ \le\ \bC Nk\Bigl((1-p)^{k^2/(2\beta)-k/2}
+(1-p)^{(k-\alpha)(\beta-\gamma)}\Bigr);
\end{equation}
if $\beta\ge\alpha$ (the condition $\beta>\gamma$ is not required in this case), then
the bound \eqref{eq:core} holds without the second term in the parentheses.
\end{theorem}

The first term in \eqref{eq:core} corresponds to Tur\'an's theorem, and the second to
Pyaderkin's lemma (Lemma~\ref{lem:pyad} in Section~\ref{sec:proofs}); the parameter
$\beta$ is chosen depending on the problem. Theorem~\ref{th:core}
is proved in Subsection~\ref{pf:core}.

\begin{corollary}
\label{sled:quant}
Let $p\in(0,1)$, $L:=\ln(eN)$, and
\begin{equation}
\label{eq:anyp}
p\ \ge\ 32L\cdot\max\Bigl\{\frac{\gamma+1}{\alpha},\ \frac1{\sqrt\alpha}\Bigr\}.
\end{equation}
Then
\begin{displaymath}
\PP\Bigl(\alpha(G_p)\ge\alpha+\frac{128L^2}{p^2}\Bigr)\le2e^{-\alpha}.
\end{displaymath}
\end{corollary}

\begin{corollary}
\label{sled:anyp}
If $(\gamma+1)\ln N=o(p\alpha)$ and $\ln^2N=o(p^2\alpha)$, then a.a.s.
\begin{displaymath}
\alpha\ \le\ \alpha(G_p)\ \le\ \alpha+\frac{128\ln^2(eN)}{p^2}=(1+o(1))\,\alpha,
\end{displaymath}
i.e., asymptotic stability holds.
\end{corollary}

\begin{corollary}
\label{sled:lambda}
Let $p\in(0,1)$, $N\to\infty$, $\lambda=\lambda_n\ge1$, $\lambda\alpha\to\infty$, and
suppose that there exists a function $\beta=\beta_n>\gamma$ such that
\begin{displaymath}
\ln N=o\Bigl(\frac{p\lambda\alpha}{\beta}\Bigr),\qquad
\ln N=o\bigl(p(\beta-\gamma)\bigr).
\end{displaymath}
Then a.a.s.\ $\alpha(G_p)\le(1+o(1))\lambda\alpha$.
\end{corollary}

\begin{corollary}
\label{sled:uncond}
Let $p\in(0,1)$, $N\to\infty$, and $\ln\frac1{1-p}=o(\ln N)$. Then a.a.s.
\begin{displaymath}
\alpha(G_p)\ \le\ (1+o(1))\cdot2\alpha\log_{1/(1-p)}N .
\end{displaymath}
If, in addition, $N/\alpha\to\infty$ and $\ln\frac1{1-p}=o(\ln(N/\alpha))$, then
a.a.s.\ $\alpha(G_p)\le(1+o(1))\cdot2\alpha\log_{1/(1-p)}(N/\alpha)$.
\end{corollary}

Corollaries~\ref{sled:quant}--\ref{sled:uncond} are derived from
Theorem~\ref{th:core} in Subsection~\ref{pf:core}.

\begin{remark}
\label{zam:A}
Theorem~A is a special case of Corollary~\ref{sled:anyp} for $p=1/2$.
Indeed, its conditions 2 and 3 imply the relations $\beta_n\ln N=o(\alpha)$ and
$\ln^2N=o(\beta_n)\cdot o(\alpha/\beta_n)=o(\alpha)$. By condition~3,
we have $\beta_n\to\infty$; therefore, $\gamma+1<\beta_n+1\le2\beta_n$ for large $n$, and
$(\gamma+1)\ln N=o(\alpha)$. Moreover, Corollary~\ref{sled:anyp} requires neither
the monotonicity of the quantities $N,\alpha,\gamma$ nor that $\gamma$ tend to
infinity; in particular, it applies to the example from \cite{Rai17} with a Steiner network, where
$\gamma\le9$ (see Subsection~\ref{ex:steiner}). In the same way, Theorem~3 of
\cite{Rai26} --- the bound $\alpha(G_{n,1/2})\le(1+o(1))\alpha_n\lambda_n$ under
the conditions of Theorem~A in which condition~2 is replaced by
$\log_2N_n=o(\alpha_n\lambda_n/\beta_n)$ --- is a special case of
Corollary~\ref{sled:lambda} for $p=1/2$, and neither the condition
$\beta_n=o(\alpha_n)$, nor monotonicity, nor $\gamma_n\to\infty$ is needed. Finally, for
$p=1/2$ the first bound of Corollary~\ref{sled:uncond} is Theorem~1 of
\cite{Rai26}, and the second is its strengthening mentioned there in Subsection~4.2 (for the graph
$G(n,3,1)$ it replaces the constant $6$ by $4$).
\end{remark}

\subsection{A Tur\'an-type bound and the order of magnitude of \texorpdfstring{$\alpha(G_p)$}{alpha(G\_p)}}

By Tur\'an's theorem \cite{Turan}, every set $D$ of cardinality $M$ contains at
least $M^2/(2\alpha)-M/2$ edges of the graph $G$. The following bound is much
stronger for large $\alpha$.

\begin{theorem}
\label{th:turan}
Let $D\subseteq V$, $M=|D|\ge16$, and $\gamma<M/2$. If $\alpha\le M/16$,
or if both $\psi\le M/16$ and $\stdeg D\le M/16$, then
$e(D)\ge M^2/16$.
\end{theorem}

Theorem~\ref{th:turan} is proved at the end of Subsection~\ref{pf:aux} (by means of
Lemma~\ref{lem:bad}).

\begin{theorem}
\label{th:order}
Let $G$ be an arbitrary graph on $N$ vertices and $p\in[0,1]$.
\begin{itemize}
\item[\rm(a)] If $k$ is a positive integer and $N(1-p)^{k-1}\ge2k$, then
$\PP(\alpha(G_p)<k)\le\exp\bigl\{-N(1-p)^{k-1}/8\bigr\}\le e^{-k/4}$.
\item[\rm(b)] If $M$ is a positive integer,
$M\ge\max\{16,\,16\alpha,\,2\gamma+1\}$, and $pM\ge16\bigl(2+\ln\frac NM\bigr)$,
then $\PP(\alpha(G_p)\ge M)\le e^{-M}$.
\end{itemize}
\end{theorem}

Bound (a) is obtained by analyzing the greedy algorithm; for the random graph
$G(n,p)$ this approach goes back to \cite{GM}, and for random subgraphs of graphs whose vertex degrees do not exceed
a sufficiently rapidly growing function $d=d(n)$, the greedy algorithm was used to
bound the chromatic number in \cite[Theorem~14]{BGPR2} and
\cite[Theorem~6]{Pyad19}; the lower bound for the independence number that follows from this for $p=1/2$
is given in \cite[Theorem~2]{Rai26}. Theorem~\ref{th:order}
is proved in Subsection~\ref{pf:order}.

\begin{corollary}
\label{sled:order}
Let $N/\alpha\to\infty$, $p\to0$, $\gamma=O(\alpha)$, and
\begin{displaymath}
f:=\frac{\ln(N/\alpha)}{p\,\alpha}\ \ge\ 1,\qquad \ln f=o\Bigl(\ln\frac N\alpha\Bigr).
\end{displaymath}
Then there exists a constant $C$, depending only on the constant in the relation
$\gamma=O(\alpha)$, such that for every fixed $\eta>0$ a.a.s.
\begin{displaymath}
(1-\eta)\,f\alpha\ \le\ \alpha(G_p)\ \le\ C f\alpha .
\end{displaymath}
\end{corollary}

Corollary~\ref{sled:order} is proved in Subsection~\ref{pf:order}.

\begin{remark}
\label{zam:Kn}
For the complete graph $K_N$ we have $\alpha=1$, $\gamma=0$, and
Corollary~\ref{sled:order} gives, for $p\to0$, $\ln(1/p)=o(\ln N)$, the well-known
order: a.a.s.\ $\alpha(G(N,p))\asymp\ln N/p$ (cf.\ the asymptotics in the introduction). If $\alpha\to\infty$, $N/\alpha\to\infty$, and $\ln(N/\alpha)=o(\alpha)$, then
Corollary~\ref{sled:order}
(with $f=1/c$) and Remark~\ref{zam:monot} show that for
$p\le c\ln(N/\alpha)/\alpha$ with a constant $c<1$ there is no asymptotic
stability. We emphasize that the lower bound in Corollary~\ref{sled:order} does not require
the condition $\gamma=O(\alpha)$.
\end{remark}

\begin{remark}
\label{zam:lower}
Bound (a) of Theorem~\ref{th:order} applies to an arbitrary graph and gives
an analogue of Theorem~4 of \cite{Rai26} for arbitrary $p$. Let $p\in(0,1)$,
let $G$ contain $f=f_n$ vertex-disjoint subgraphs, each
of which has at least $w=w_n$ vertices, with no edges between different
subgraphs, and let $w\to\infty$,
$\ln\frac1p=o(\ln w)$, $\ln\frac1{1-p}=o(\ln w)$, and $\ln f=O(w^{\varepsilon})$
for every fixed $\varepsilon>0$. Then for every fixed
$\varepsilon\in(0,1)$ a.a.s.
\begin{displaymath}
\alpha(G_p)\ \ge\ (1-\varepsilon)\,f\log_{1/(1-p)}w .
\end{displaymath}
Indeed, set $k:=\lfloor(1-\varepsilon/2)\log_{1/(1-p)}w\rfloor$;
then $w(1-p)^{k-1}\ge w^{\varepsilon/2}\ge2k$ for large $n$ (since
$k\le\ln w/p\le w^{\varepsilon/4}$), and by (a) the probability that the random
subgraph of at least one of the $f$ subgraphs contains no independent set
of cardinality $k$ does not exceed $f\exp\{-w^{\varepsilon/2}/8\}\to0$. The union of
independent sets taken from distinct subgraphs, one from each, is independent in
$G_p$; therefore, a.a.s.\ $\alpha(G_p)\ge fk\ge(1-\varepsilon)f\log_{1/(1-p)}w$
(the last inequality holds for large $n$, since
$\log_{1/(1-p)}w\to\infty$). For $p=1/2$, in \cite[Theorem~4]{Rai26}, owing to the
exact asymptotics of the independence number of the Erd\H{o}s--R\'enyi random graph, a
bound better by a factor of two, $(1+o(1))\cdot2f\log_2w$ is obtained (under the conditions $f\to\infty$,
$w\to\infty$, $fe^{-w}\to0$). The lower bound in terms of the maximum vertex degree (Theorem~2
of \cite{Rai26} for $p=1/2$) follows for arbitrary $p$ from
\cite[Theorem~6]{Pyad19}.
\end{remark}

\subsection{Threshold for asymptotic stability}

\begin{theorem}
\label{th:asthr}
Let $\alpha\to\infty$, $N/\alpha\to\infty$, $\ln(N/\alpha)=o(\alpha)$, and
suppose that for all sufficiently large $n$
\begin{equation}
\label{eq:struct}
\gamma<\frac\alpha2,\qquad \psi\le\frac{\alpha}{16},\qquad
\delta\le\frac1{32}.
\end{equation}
Set $p^*:=\ln(N/\alpha)/\alpha$. Then the following
statements hold.
\begin{itemize}
\item[\rm(i)] For every fixed $\varepsilon\in(0,1]$, if
$p\ge192\ln(eN/\alpha)/(\varepsilon\alpha)$, then a.a.s.\
$\alpha(G_p)\le(1+\varepsilon)\alpha$. In particular, if $p/p^*\to\infty$, then
a.a.s.\ $\alpha(G_p)=(1+o(1))\alpha$.
\item[\rm(ii)] For any fixed $c\in(0,1)$ and $\eta>0$, if $p\le cp^*$, then
a.a.s.\ $\alpha(G_p)\ge(1-\eta)\alpha/c$; this statement
also holds without assumptions \eqref{eq:struct}.
\end{itemize}
\end{theorem}

Thus, $p^*$ is the threshold probability for asymptotic
stability: for $p\gg p^*$ it holds, for $p\le cp^*$, $c<1$, it
fails, and for
$p\asymp p^*$ the independence number $\alpha(G_p)$ is a.a.s.\ of order $\alpha$ (the upper
bound $O(\alpha)$ follows from Corollary~\ref{sled:order} and
Remark~\ref{zam:monot}). Theorem~\ref{th:asthr} is proved in
Subsection~\ref{pf:asthr}.

\subsection{Threshold for exact stability}

\begin{theorem}
\label{th:exthr}
Let $\alpha\to\infty$ and $\Lambda\to\infty$.
\begin{itemize}
\item[\rm(i)] Suppose that for all sufficiently large $n$
\begin{equation}
\label{eq:struct2}
\gamma<\frac\alpha2,\qquad \psi\le\frac{\alpha}{16},\qquad
\delta\le\frac1{64}.
\end{equation}
Then for $p\ge p_1:=64\Lambda/\alpha$ a.a.s.\ $\alpha(G_p)=\alpha$ and
$\I(G_p)=\I(G)$, i.e., the independent sets of maximum cardinality of the graph $G_p$
are exactly the independent sets of maximum cardinality of the graph $G$.
\item[\rm(ii)] Let $\Lambda^2=o(\alpha)$ and
\begin{equation}
\label{eq:second}
\ln(N-\alpha)-4\Lambda\Bigl(\delta+\frac\theta\alpha\Bigr)\to+\infty .
\end{equation}
Then for every fixed $\varepsilon\in(0,1)$ and every
$p\le p_2:=(1-\varepsilon)\Lambda/\alpha$, a.a.s.\
$\alpha(G_p)\ge\alpha+1$.
\end{itemize}
\end{theorem}

Part (i) of Theorem~\ref{th:exthr} is proved in Subsection~\ref{pf:exthr1}, and part
(ii) in Subsection~\ref{pf:exthr2}.

\begin{corollary}
\label{sled:exthr}
If the conditions of both parts of Theorem~\ref{th:exthr} are satisfied, then the threshold
probability for exact stability has order $\Lambda/\alpha$, and the ratio
of the upper bound to the lower bound equals $p_1/p_2=64/(1-\varepsilon)$.
\end{corollary}

Corollary~\ref{sled:exthr} follows directly from
Theorem~\ref{th:exthr}.

\begin{remark}
\label{zam:compare}
Let us compare the thresholds. Since $\ln(N-\alpha)\le\Lambda$, and for $\alpha\le N/2$ and
$\alpha\ge2$ also $\ln(N/\alpha)\le\ln(N/2)\le\ln(N-\alpha)$, the threshold
for asymptotic stability $p^*$ does not exceed the threshold for exact stability
$\Lambda/\alpha$. The ratio $\Lambda/\ln(N/\alpha)$ can be either
bounded (Subsection~\ref{ex:steiner}) or unbounded (for example, in
Subsection~\ref{ex:multi} for $t=\lfloor\ln b\rfloor$, in Subsection~\ref{ex:Gnrs} for $s=1$,
$r\to\infty$, in Subsection~\ref{ex:qkneser} for $k\to\infty$); in the second case, if
the conditions of Theorems~\ref{th:asthr} and~\ref{th:exthr}(ii) are satisfied, then for
$p^*\ll p\le(1-\varepsilon)\Lambda/\alpha$ a.a.s.\ $\alpha(G_p)=(1+o(1))\alpha$,
but $\alpha(G_p)>\alpha$.
\end{remark}

\section{Examples}
\label{sec:ex}

Throughout this section, $\varepsilon\in(0,1)$ is an arbitrary fixed
number.

\subsection{Complete multipartite graphs}
\label{ex:multi}

Let $G=K_{t\times b}$ be the complete $t$-partite graph with parts of cardinality $b$,
$t\ge2$, $b\ge2$. The independent sets are the subsets of the parts; therefore,
$\alpha=b$, $\I$ is the set of parts, $m=t$, $\psi=0$, $\gamma=0$ (a common
non-neighbor of the endpoints of an edge must lie in the part of each of the endpoints), $\mu=b$,
$\delta=0$, $\theta=0$, $N=tb$, $\Lambda=\ln\bigl(t(t-1)b\bigr)$.

\begin{proposition}
\label{utv:multi}
Let $b\to\infty$ and $\ln t=o(\sqrt b)$. Then for $p\ge64\ln(t^2b)/b$
a.a.s.\ $\alpha(G_p)=b$ and the independent sets of cardinality $b$ in the graph $G_p$
are the parts, while for $p\le(1-\varepsilon)\ln(t(t-1)b)/b$ a.a.s.\ $\alpha(G_p)>b$.
If, in addition, $t\to\infty$, then the threshold for asymptotic stability
(in the sense of Theorem~\ref{th:asthr}) equals $p^*=\ln t/b$.
\end{proposition}

\begin{proof}
Conditions \eqref{eq:struct2} are satisfied, $\Lambda\to\infty$, $\Lambda^2\le(2\ln t+\ln b)^2=o(b)$, and \eqref{eq:second} reduces to
$\ln((t-1)b)\to\infty$. It remains to apply Theorems~\ref{th:exthr}
and~\ref{th:asthr}, noting that $\ln(t(t-1)b)\le\ln(t^2b)$.
\end{proof}

For $t=\lfloor\ln b\rfloor$ we have $p^*\asymp\ln\ln b/b$, whereas
$\Lambda/\alpha\asymp\ln b/b$. For comparison, Corollary~\ref{sled:anyp} gives
asymptotic stability here only for $p\gg\ln(tb)/\sqrt b$.

\subsection{The construction from~\cite{Rai17}}
\label{ex:bip}

Let $V=A\sqcup B$, $|A|=a\ge1$, $|B|=b>a$, where the set $B$ is independent, each
vertex of $A$ is adjacent to each vertex of $B$, and the edges inside $A$
are arbitrary. An independent set meeting $A$ contains no vertices of
$B$, i.e., lies in $A$. Therefore, $\alpha=b$, $\I=\{B\}$, $m=1$, $\psi\le a$,
$\theta=0$, $\mu=b$, $\delta=0$, $N-\alpha=a$, $\Lambda=\ln a$. Further,
$\gamma\le a$: for an edge inside $A$ the common non-neighbors lie in $A$, and an edge
$\{x,y\}$ with $x\in A$, $y\in B$ has no common non-neighbors.

\begin{proposition}
\label{utv:bip}
Let $a\to\infty$ and $16a\le b$. Then for $p\ge64\ln a/b$ a.a.s.\
$\alpha(G_p)=b$ and $B$ is the unique independent set of cardinality $b$ in
$G_p$, while for $p\le(1-\varepsilon)\ln a/b$ a.a.s.\ $\alpha(G_p)>b$.
\end{proposition}

\begin{proof}
This is Theorem~\ref{th:exthr}: conditions \eqref{eq:struct2} are satisfied, since
$\gamma\le a\le b/16$, $\psi\le a\le b/16$, $\delta=0$; further,
$\Lambda^2=\ln^2a=o(a)=o(b)$, and \eqref{eq:second} means $\ln a\to\infty$.
\end{proof}

The lower bound here is also easy to obtain directly: the probability that
at least one vertex of $A$ is joined in $G_p$ to no vertex of $B$
(and then $\alpha(G_p)>b$) equals $1-\bigl(1-(1-p)^b\bigr)^a$, and this tends
to one for $p\le(1-\varepsilon)\ln a/b$, $p\to0$. Proposition~\ref{utv:bip}
shows that the order $\ln a/b$ is sharp. In \cite{Rai17}, for this construction it was assumed that
$a=o(\sqrt b)$, and asymptotic stability was considered for $p=1/2$.
The same construction is analyzed in \cite[Subsection~4.1]{Rai26}: for $\lambda_n=1$, the
condition $a=o(b/\log_2b)$ is obtained there; Corollary~\ref{sled:anyp} for $p=1/2$ gives
the condition $(\gamma+1)\ln b=o(b)$, which follows from $a=o(b/\log_2b)$, since
$\gamma\le a$, and coincides with it if $\gamma\asymp a$ (nothing is assumed about the edges inside $A$
in \cite{Rai26}); for $\lambda_n>1$, larger values of $a$ are also allowed in \cite{Rai26}.
Corollary~\ref{sled:lambda} gives here, for arbitrary $p$,
the following: if $a\ln b=o(p\lambda b)$ and $\ln^2b=o(p^2\lambda b)$, then a.a.s.\
$\alpha(G_p)\le(1+o(1))\lambda b$ (it suffices to take $\beta=a+\omega\ln b/p$ with
a suitable $\omega\to\infty$).

\subsection{The graphs \texorpdfstring{$G(n,r,<s)$}{G(n,r,<s)}}
\label{ex:Gnrs}

Let $1\le s<r<n$, let $V=V(n,r)$ be the family of all $r$-element
subsets of the set $[n]$, and let vertices $x,y$ be adjacent if $|x\cap y|<s$.
We set (unlike \cite{PR2}, where the independence number is denoted
by $N$)
\begin{displaymath}
N=\bC nr,\qquad a:=\bC{n-s}{r-s},\qquad c:=\bC rs .
\end{displaymath}
The \emph{star} $S_W=\{x\in V\colon x\supseteq W\}$ centered at $W\subset[n]$,
$|W|=s$, is independent and has cardinality $a$. Stars with different centers are
different (the center is the intersection of all elements of the star).

\begin{proposition}
\label{utv:Gnrs}
Let $r\ge2$, and suppose that for all sufficiently large $n$ at least one of the
conditions
\begin{displaymath}
\text{\rm(A)}\ \ n\ge32\,c\,r^2,\qquad\text{\rm(B)}\ \ n\ge30\,s\,r^3
\end{displaymath}
holds (this is the case, in particular, if $\min\{cr^2,sr^3\}=o(n)$). Then $\alpha=a$,
$\I$ is the family of all stars, $m=\bC ns$,
$\Lambda=\ln\bigl(\bC ns(N-a)\bigr)$, and the following holds.
\begin{itemize}
\item[\rm(a)] For $p\ge64\Lambda/a$ a.a.s.\ $\alpha(G_p)=a$ and all
independent sets of cardinality $a$ in $G_p$ are stars; for
$p\le(1-\varepsilon)\Lambda/a$ a.a.s.\ $\alpha(G_p)\ge a+1$.
\item[\rm(b)] The threshold for asymptotic stability (in the sense of
Theorem~\ref{th:asthr}) equals $p^*=\ln\bigl(\bC ns/\bC rs\bigr)/\bC{n-s}{r-s}$.
\end{itemize}
\end{proposition}

\begin{proof}
We use the bounds from \cite{PR2}, where the quantity $\gamma$ and
the ``nontrivial'' analogue of the quantity $\psi$ are denoted by $\Phi$ and $\Psi$:
there, $\Phi$ is the maximum number of common non-neighbors of the endpoints of an edge, and $\Psi$ is
the maximum cardinality of an $s$-intersecting family not contained in any
star. In the proof of Corollary~2 of \cite{PR2}, it is shown that in
each of the cases (A) and (B) one has $\Phi<a/2$ and $\Psi\le a/16$ (in
case (B), the Ahlswede--Khachatrian theorem \cite{AK2} is used for this).
An independent set either lies in a star and then has cardinality at most
$a$, or has cardinality at most $\Psi<a$. Hence $\alpha=a$, $\I$ consists
exactly of the stars, $m=\bC ns$, $\psi=\Psi\le a/16$, $\gamma=\Phi<a/2$.
Further, by Lemma~2 of \cite{PR2}, for large $n$ every vertex outside a star
is adjacent to at least $(1-2r^2/n)a$ vertices of this star, i.e.,
$\delta\le2r^2/n$. Finally, two distinct stars intersect in at most
$\bC{n-s-1}{r-s-1}=a\frac{r-s}{n-s}\le\frac{ar}{n}$ vertices, so that
$\theta/a\le r/n$.

In each of the cases (A) and (B) one has $n\ge30r^3$ (in case (A), since
$c\ge r$ for $1\le s\le r-1$); therefore, $r^2\le n^{2/3}$. Hence
$\delta\le2r^2/n\to0$, and conditions \eqref{eq:struct2} are satisfied; moreover,
$r^2\ln^2n=o(n)$. Since $a/N=c/\bC ns\le(er/n)^s\to0$, for large $n$
we have $N-a\ge N/2$. Using the bounds $(x/y)^y\le\bC xy\le(ex/y)^y$, we obtain
\begin{displaymath}
L-\ln2\ \le\ \Lambda\ \le\ L+r+s\ \le\ 2r(\ln n+1),\qquad
L:=r\ln\frac nr+s\ln\frac ns .
\end{displaymath}
Since $L\ge r\ln\frac nr\ge\frac r2\ln n\to\infty$, we have
$\Lambda=(1+o(1))L\to\infty$. Further, $a\ge n-s\ge n-r$; therefore,
$\Lambda^2\le4r^2(\ln n+1)^2=o(a)$. Finally,
$\ln(N-a)\ge\ln(N/2)\ge\frac r2\ln n-\ln2$, whereas
\begin{displaymath}
4\Lambda\Bigl(\delta+\frac\theta a\Bigr)\le8r(\ln n+1)\cdot\frac{3r^2}{n}
=o(r\ln n),
\end{displaymath}
i.e., \eqref{eq:second} holds. Statement (a) follows from
Theorem~\ref{th:exthr}. For (b), note that $N/a=\bC ns/c$ (the identity
$\bC nr\bC rs=\bC ns\bC{n-s}{r-s}$), and moreover $\bC ns/c\ge(n/r)^s\to\infty$ and
$\ln(N/a)\le s\ln(en/s)\le r\ln(en)=o(a)$; conditions \eqref{eq:struct} follow
from \eqref{eq:struct2}, and it remains to apply Theorem~\ref{th:asthr}.
\end{proof}

Proposition~\ref{utv:Gnrs}(a) contains Theorem~1 of \cite{PR2}
(there the lower bound is written as $(1-\varepsilon)L/a$; since
$\Lambda=(1+o(1))L$, both forms are equivalent up to a change of
$\varepsilon$). Proposition~\ref{utv:Gnrs}(b) is apparently new.

\emph{Kneser graphs} ($s=1$, condition (A) takes the form $n\ge32r^3$).
Here $p^*=\ln(n/r)\big/\bC{n-1}{r-1}$, whereas the threshold for exact
stability found in \cite{BNR} for $2\le r=o(n^{1/3})$ up to a
factor $1+o(1)$ equals $p_c:=\bigl((r+1)\ln n-r\ln r\bigr)\big/\bC{n-1}{r-1}$.
The ratio $p_c/p^*$ is greater than $r$; therefore, if $r\to\infty$ (in the range
$2\le r=o(n^{1/3})$, where both results apply), then for
$p^*\ll p\le(1-\varepsilon)p_c$ a.a.s.\ $\alpha(G_p)\sim\bC{n-1}{r-1}$, but
$\alpha(G_p)>\bC{n-1}{r-1}$.

\subsection{Steiner systems}
\label{ex:steiner}

Let $k\ge2$, and let $\mathcal S$ be a Steiner system $S(2,k,v)$, i.e., a family
of $k$-element subsets (\emph{blocks}) of a $v$-element set of points
in which each pair of points lies in exactly one block. Through each point
there pass $r=(v-1)/(k-1)$ blocks, and the total number of blocks is $N=v(v-1)/(k(k-1))$.
Consider the graph whose vertices are the blocks and whose edges are the pairs
of disjoint blocks; for $k=3$ (Steiner triple systems) this is the graph from the
example of \cite{Rai17}, where such a system is called a Steiner network. The \emph{star} $S_x$ is the set of blocks containing the point
$x$; it is independent and has cardinality $r$.

Systems $S(2,3,v)$ exist for all $v\equiv1,3\pmod 6$ \cite{Kirkman}
(see also the books \cite{Hall,ColbournRosa}); for every fixed $k$,
systems $S(2,k,v)$ exist for all sufficiently large $v$
satisfying the necessary divisibility conditions \cite{Wilson75}; the lines
of the affine space $AG(d,q)$ form a system $S(2,q,q^d)$. An analogous
existence theorem for Steiner systems $S(t,k,v)$ and, more generally, for $t$-designs
with arbitrary fixed $t<k$ was proved relatively
recently by P.~Keevash \cite{Keevash} (for another proof, see \cite{GKLO}); for a survey of known constructions, see the handbook
\cite{CD}.

\begin{proposition}
\label{utv:steiner}
Let $v\to\infty$, $k=k(v)\ge2$, and $64k^2(k-1)\le v-1$, i.e., $64k^2\le r$.
Then $\alpha=r$, $\I$ is the family of stars, $m=v$, and
\begin{displaymath}
\gamma\le k^2,\qquad \psi\le k^2-k+1,\qquad \delta\le\frac kr,\qquad
\theta=1,\qquad 2\ln v\le\Lambda\le3\ln v .
\end{displaymath}
Consequently, for $p\ge64\Lambda/r$ (in particular, for
$p\ge192(k-1)\ln v/(v-1)$) a.a.s.\ $\alpha(G_p)=r$ and all independent
sets of cardinality $r$ in $G_p$ are stars; for $p\le(1-\varepsilon)\Lambda/r$
a.a.s.\ $\alpha(G_p)\ge r+1$; finally, for
$p\le c\ln(v/k)\,(k-1)/(v-1)$, $c\in(0,1)$, for every $\eta>0$ a.a.s.\
$\alpha(G_p)\ge(1-\eta)r/c$. Both thresholds have order $(k-1)\ln v/v$.
\end{proposition}

\begin{proof}
Two distinct blocks have at most one common point. Let $T,T'$ be
disjoint blocks, and let a block $T''$ meet both; take $u\in T\cap
T''$, $w\in T'\cap T''$; then $T''$ is the unique block containing $u$ and
$w$. Therefore, $\gamma\le k^2$. If a block $T$ does not contain a point $x$, then
the number of blocks through $x$ meeting $T$ is at most $k$ (one through $x$ and each
point of $T$); therefore, $T$ is adjacent to at least $r-k$ blocks of the star $S_x$.
The block containing two distinct points $x,y$ is unique, so that $|S_x\cap
S_y|=1$.

Let $\mathcal B$ be a family of pairwise intersecting blocks not contained
in any star. Take $T_1,T_2\in\mathcal B$, $T_1\cap T_2=\{x\}$, and a
block $T_3\in\mathcal B$ not containing $x$. Every block $T\in\mathcal B$
containing $x$ meets $T_3$ in a point different from $x$ and is uniquely
determined by this point --- at most $k$ possibilities. Every block
$T\in\mathcal B$ not containing $x$ meets $T_1$ and $T_2$ in distinct
points $u\in T_1\setminus\{x\}$, $w\in T_2\setminus\{x\}$ and is uniquely
determined by the pair $(u,w)$ --- at most $(k-1)^2$ possibilities. Thus,
$|\mathcal B|\le k+(k-1)^2=k^2-k+1<r$. As in Subsection~\ref{ex:Gnrs}, this implies
$\alpha=r$, $\I$ is the family of stars (the stars are distinct, since $r\ge2$),
$\psi\le k^2-k+1$, $\mu\ge r-k$.

Further, $N\le v^2/2$ and $N-r\ge N/2$ (since $N/r=v/k\ge2$); therefore,
$\ln(vN/2)\le\Lambda\le\ln(v^3/2)$, and $vN/2=v^2(v-1)/(2k(k-1))\ge v^2$
since $64k^2\le v-1$; this gives $2\ln v\le\Lambda\le3\ln v$.

Conditions \eqref{eq:struct2}: $\gamma\le k^2\le r/64$, $\psi\le k^2\le r/64$,
$\delta\le k/r\le k^2/r\le1/64$. From $64k^2(k-1)\le v-1$ it follows that
$k-1\le(v-1)^{1/3}/4$ and $r\ge4(v-1)^{2/3}\to\infty$; therefore,
$\Lambda^2\le9\ln^2v=o(r)$; moreover, $k^2\le(v-1)^{2/3}$ for large $v$;
therefore,
\begin{displaymath}
4\Lambda\Bigl(\delta+\frac1r\Bigr)\le12\ln v\cdot\frac{k+1}r\le
12\ln v\cdot\frac{k^2}{v-1}\le\frac{12\ln v}{(v-1)^{1/3}}\to0,
\end{displaymath}
and $\ln(N-r)\to\infty$, i.e., \eqref{eq:second} holds. The statements on exact stability follow from
Theorem~\ref{th:exthr}. Finally, $N/\alpha=v/k\to\infty$ and
$\ln(v/k)=o(r)$; therefore, Theorem~\ref{th:asthr} applies with
$p^*=\ln(v/k)/r=\ln(v/k)(k-1)/(v-1)$; since $k\le v^{1/3}$, we have
$\ln(v/k)\ge\frac23\ln v$, and both thresholds have order $(k-1)\ln v/v$.
\end{proof}

For a Steiner triple system ($k=3$), this graph is given in \cite{Rai17} as an example
of asymptotic stability for $p=1/2$ (formally, it does not satisfy the condition
$\gamma_n\to\infty$ of Theorem~A, see Remark~\ref{zam:A}); Proposition~\ref{utv:steiner} gives exact
stability already for $p\ge384\ln v/(v-1)$, and the order of this threshold cannot be
improved.

\subsection{\texorpdfstring{$q$}{q}-analogues of Kneser graphs}
\label{ex:qkneser}

Let $q$ be a prime power, $2\le k<n/2$, and let $V$ be the set of all
$k$-dimensional subspaces of the space $\mathbb F_q^n$; subspaces are adjacent
if their intersection is zero. Denote by $\qb ab$ the Gaussian binomial
coefficients (the number of $b$-dimensional subspaces of an $a$-dimensional space over
$\mathbb F_q$); the number of $k$-dimensional subspaces containing a given
$j$-dimensional one equals $\qb{n-j}{k-j}$. The \emph{star} $S_P$ is the set of
subspaces containing a given one-dimensional subspace (\emph{point})
$P$; it is independent and has cardinality $a:=\qb{n-1}{k-1}$. Set
$\varkappa:=\qb k1^2\qb{n-2}{k-2}$.

\begin{proposition}
\label{utv:q}
Let $n\to\infty$ (the quantities $q\ge2$ and $k\ge2$ may depend on $n$) and
$n\ge3k+5$. Then $\alpha=a$, $\I$ is the family of stars, $m=\qb n1$, and
\begin{displaymath}
\gamma\le\varkappa\le\frac a{32},\qquad\psi\le1+\varkappa\le\frac a{16},\qquad
\delta\le\frac1{128},\qquad \theta=\qb{n-2}{k-2}\le q^{k-n}a .
\end{displaymath}
Consequently, the threshold for exact stability has order $\Lambda/a$, where
$\Lambda=\ln\bigl(\qb n1(\qb nk-a)\bigr)\asymp kn\ln q$, and the threshold for asymptotic
stability equals $p^*=\ln\frac{q^n-1}{q^k-1}\big/a\asymp n\ln q/a$.
\end{proposition}

\begin{proof}
\emph{Bound on $\gamma$.} Let $U\cap U'=\{0\}$, and let a subspace $W$
intersect nontrivially both $U$ and $U'$. Then $W$ contains a point $P\subseteq
U$ and a point $P'\subseteq U'$, with $P\ne P'$, i.e., it contains the two-dimensional
subspace $P+P'$. The number of pairs $(P,P')$ is at most $\qb k1^2$, and for each pair
the number of subspaces $W\supseteq P+P'$ is exactly $\qb{n-2}{k-2}$. Hence
$\gamma\le\varkappa$.

\emph{Bound on $\psi$.} Let $\mathcal B$ be a family of pairwise nontrivially
intersecting subspaces not contained in any star, and let
$X\in\mathcal B$. Every $Y\in\mathcal B\setminus\{X\}$ contains some
point $P\subseteq X$ (at most $\qb k1$ possibilities). For fixed $P$
there exists $Z_P\in\mathcal B$ not containing $P$; then every $Y\in\mathcal B$
containing $P$ is different from $Z_P$ and contains a point $P'\subseteq Z_P$, $P'\ne
P$, and hence also the subspace $P+P'$; the number of such $Y$ is at most
$\qb k1\qb{n-2}{k-2}$. Thus, $|\mathcal B|\le1+\varkappa$.

\emph{Bound on $\delta$.} If $P\not\subseteq U$, then the subspaces $W\in
S_P$ intersecting $U$ nontrivially contain $P+Q$ for some point
$Q\subseteq U$; there are at most $\qb k1\qb{n-2}{k-2}$ of them. Therefore,
$\delta\le\qb k1\qb{n-2}{k-2}/a$.

\emph{Numerical bounds.} We have $a=\qb{n-2}{k-2}\frac{q^{n-1}-1}{q^{k-1}-1}$,
$\frac{q^{k-1}-1}{q^{n-1}-1}\le q^{k-n}$, and $\qb k1=\frac{q^k-1}{q-1}\le
2q^{k-1}$. Therefore, using the inequalities $n\ge3k+5$ and $q\ge2$, we obtain
\begin{displaymath}
\frac\varkappa a\le4q^{3k-2-n}\le4q^{-7}\le\frac1{32},\qquad
\delta\le2q^{2k-1-n}\le2q^{-k-6}\le\frac1{128},\qquad
\frac\theta a\le q^{k-n}.
\end{displaymath}
Moreover, $a\ge\qb{n-1}1\ge2^{n-2}\ge32$, so that $1+\varkappa\le a/16$.
Since $1+\varkappa<a$, an independent set either lies in a star or
has cardinality less than $a$; hence $\alpha=a$, $\I$ is the family of
(distinct) stars, $m=\qb n1$, and $\psi\le1+\varkappa$. The quantity $|S_P\cap S_Q|$ for
$P\ne Q$ equals $\qb{n-2}{k-2}$.

\emph{Conditions of the theorems.} Conditions \eqref{eq:struct2} (and hence also
\eqref{eq:struct}) have already been verified. Further, $m\le q^n$, $\qb nk\le q^{nk}$
(a subspace is determined by a basis), and $\qb nk-a\ge\frac12\qb nk\ge\frac12
q^{k(n-k)}$; therefore, $k(n-k)\ln q-\ln2\le\Lambda\le(k+1)n\ln q$, whence
$\Lambda\asymp kn\ln q$ (since $n-k\ge 2n/3$). Set $y:=(k-1)(n-k)\ln q$;
then $a\ge q^{(k-1)(n-k)}=e^y$, $y\to\infty$, and
$\Lambda\le(k+1)n\ln q\le\frac92y$, so that $\Lambda^2=o(a)$. For
\eqref{eq:second}, set $u:=n-2k+1$; then $u\ge k+6$, $u\ge n/3$, and
\begin{displaymath}
4\Lambda\Bigl(\delta+\frac\theta a\Bigr)\le4(k+1)n\ln q\cdot3q^{-u}
\le36\,u^2\ln q\cdot q^{-u}\to0,
\end{displaymath}
and $\ln(\qb nk-a)\to\infty$. Finally, $N/a=\frac{q^n-1}{q^k-1}\ge q^{n-k}\to\infty$
and $\ln(N/a)\le n\ln q=o(a)$. It remains to apply Theorems~\ref{th:asthr}
and~\ref{th:exthr}.
\end{proof}

The equality $\alpha=\qb{n-1}{k-1}$ and the fact that the stars are the only independent
sets of maximum cardinality constitute (in the parameter range under consideration) the
Erd\H{o}s--Ko--Rado theorem for vector spaces \cite{Hsieh,FW}; here it is
obtained anew from the bound $\psi\le1+\varkappa$. For $k=2$, the graph from
Proposition~\ref{utv:q} is the graph of disjoint lines of a projective
space, i.e., the graph on the blocks of a Steiner system as in Subsection~\ref{ex:steiner}.

\section{Proofs}
\label{sec:proofs}

\subsection{Auxiliary results and proof of Theorem~\ref*{th:turan}}
\label{pf:aux}

We repeatedly use the fact that a set $D$ is independent in $G_p$ with
probability $(1-p)^{e(D)}\le e^{-pe(D)}$, the inequality
$\bC Nk\le(eN/k)^k$, and Markov's inequality $\PP(X>0)\le\EE X$ for
nonnegative integer-valued $X$.

\begin{lemma}[\cite{BGPR1,BGPR2,Rai17}]
\label{lem:pyad}
Let $A\subseteq V$, and let $B$ be a maximum-cardinality independent
subset of the set $A$. Then $e(A)\ge(|A|-|B|)(|B|-\gamma)$.
\end{lemma}

\begin{proof}
Let $x\in A\setminus B$. The set $B\cup\{x\}$ is not independent; therefore, $x$
is adjacent to some vertex $y\in B$. Every vertex $z\in B$ not adjacent to
$x$ is different from $y$ and is not adjacent to $y$ (since $B$ is independent), i.e.,
it is a common non-neighbor of the endpoints of the edge $\{x,y\}$; the number of such $z$ is at most $\gamma$.
Hence $x$ is adjacent to at least $|B|-\gamma$ vertices of $B$. The edges
joining distinct $x\in A\setminus B$ to $B$ are distinct.
\end{proof}

\begin{lemma}[Tur\'an \cite{Turan}]
\label{lem:turan}
If $A\subseteq V$ and $b=\alpha(G[A])\ge1$, then $e(A)\ge|A|^2/(2b)-|A|/2$.
\end{lemma}

\begin{lemma}
\label{lem:bad}
Let $D\subseteq V$, $M=|D|$, $\gamma<M/2$. For $x\in D$, denote by
$\nu(x)$ the number of vertices of the set $D\setminus\{x\}$ not adjacent to $x$, and
set $B:=\{x\in D\colon\nu(x)\ge\frac34M\}$. Then $B$ is independent, and if
$M\ge16$ and $|B|\le M/16$, then $e(D)\ge M^2/16$.
\end{lemma}

\begin{proof}
Suppose that $x,x'\in B$ are adjacent, and let $D_x$, $D_{x'}$ be the sets of non-neighbors
of the vertices $x$, $x'$ in $D\setminus\{x\}$ and $D\setminus\{x'\}$, respectively.
Both sets lie in $D$; therefore,
$|D_x\cap D_{x'}|\ge\nu(x)+\nu(x')-M\ge M/2$. Every vertex of $D_x\cap
D_{x'}$ is different from $x$ and $x'$ and is adjacent to neither of them, so that
$|D_x\cap D_{x'}|\le\gamma<M/2$ --- a contradiction. Hence $B$ is independent.

Further, $\nu(x)\le M-1$ for $x\in B$ and $\nu(x)<\frac34M$ for $x\notin B$.
Therefore, for $|B|\le M/16$
\begin{displaymath}
e(D)=\bC M2-\frac12\sum_{x\in D}\nu(x)\ \ge\ \frac{M(M-1)}2-\frac12\Bigl(
\frac{M^2}{16}+\frac34M^2\Bigr)=\frac{3M^2}{32}-\frac M2\ \ge\ \frac{M^2}{16},
\end{displaymath}
where the last inequality is equivalent to $M\ge16$.
\end{proof}

\begin{proof}[Proof of Theorem~\ref{th:turan}]
Let $B$ be the set from Lemma~\ref{lem:bad}; it is independent. If
$\alpha\le M/16$, then $|B|\le\alpha\le M/16$. Let $\psi\le M/16$ and
$\stdeg D\le M/16$. If $B\subseteq I$ for some $I\in\I$, then, since
$B\subseteq D$, we have $|B|\le|D\cap I|\le\stdeg D\le M/16$; otherwise
$|B|\le\psi\le M/16$ by the definition of $\psi$. In both cases, it remains to
apply Lemma~\ref{lem:bad}.
\end{proof}

\subsection{Proof of Theorem~\ref*{th:core} and Corollaries~\ref*{sled:quant}--\ref*{sled:uncond}}
\label{pf:core}

\emph{Theorem~\ref{th:core}.} If $\alpha(G_p)\ge k$, then $G_p$ contains an
independent set $A$ of cardinality $k$. Let $B$ be a maximum-cardinality subset of the set $A$
that is independent in $G$, and let $b:=|B|$; then
$1\le b\le\alpha<k$. If $b\le\beta$, then by Lemma~\ref{lem:turan} we have
$e(A)\ge k^2/(2b)-k/2\ge k^2/(2\beta)-k/2$. If $b>\beta$ (for
$\beta\ge\alpha$ this case is impossible, since $b\le\alpha$, and then
the second term is not needed), then, since $\beta>\gamma$, by Lemma~\ref{lem:pyad} we obtain
$e(A)\ge(k-b)(b-\gamma)\ge(k-\alpha)(\beta-\gamma)$, because
$k-b\ge k-\alpha>0$ and $b-\gamma>\beta-\gamma>0$. The set $A$ is independent in
$G_p$ with probability $(1-p)^{e(A)}$; summing over all $A$ of cardinality $k$,
we obtain \eqref{eq:core}. \hfill$\square$

\emph{Corollary~\ref{sled:quant}.} Set $t:=\lceil128L^2/p^2\rceil$, $k:=\alpha+t$, and
$\beta:=\lfloor p\alpha/(8L)\rfloor$. It follows from \eqref{eq:anyp} that
$p\alpha/(8L)\ge4(\gamma+1)\ge4$; therefore,
\begin{displaymath}
\beta\ \ge\ \frac{p\alpha}{8L}-1\ \ge\ \frac{p\alpha}{16L}\ \ge\ 2(\gamma+1),
\qquad\text{whence}\qquad \beta-\gamma\ \ge\ \frac\beta2\ \ge\
\frac{p\alpha}{32L}.
\end{displaymath}
Moreover, it follows from \eqref{eq:anyp} that $128L^2/p^2\le\alpha/8$; therefore, $t\le\alpha$ and
$k\le2\alpha$. We may assume that $k\le N$ (otherwise $\PP(\alpha(G_p)\ge k)=0$).
Since $\beta>\gamma$ and $(k-\alpha)(\beta-\gamma)=t(\beta-\gamma)$, by
Theorem~\ref{th:core} and the inequality $(1-p)^x\le e^{-px}$
\begin{displaymath}
\PP(\alpha(G_p)\ge k)\ \le\
\bC Nk\Bigl(e^{-p(k^2/(2\beta)-k/2)}+e^{-pt(\beta-\gamma)}\Bigr),
\end{displaymath}
and $\bC Nk\le(eN/k)^k\le e^{kL}$. Since $\beta\le p\alpha/(8L)$ and $k\ge\alpha$,
we have $pk/(2\beta)\ge4L$, whence
\begin{displaymath}
kL-p\Bigl(\frac{k^2}{2\beta}-\frac k2\Bigr)\ \le\ kL-4kL+\frac k2\ \le\ -k\ \le\
-\alpha
\end{displaymath}
(here $L\ge1$). Further, $p^2t\ge128L^2$ and $k\le2\alpha$; therefore,
\begin{displaymath}
kL-pt(\beta-\gamma)\ \le\ kL-\frac{p^2t\,\alpha}{32L}\ \le\ kL-4\alpha L\
\le\ -2\alpha L\ \le\ -\alpha .
\end{displaymath}
Thus, $\PP(\alpha(G_p)\ge k)\le2e^{-\alpha}$, and the event
$\{\alpha(G_p)\ge\alpha+128L^2/p^2\}$ coincides with $\{\alpha(G_p)\ge k\}$, since
$\alpha(G_p)$ is an integer. \hfill$\square$

\emph{Corollary~\ref{sled:anyp}.} The condition
$\ln^2N=o(p^2\alpha)$ and the inequality $N\ge2$ imply $\alpha\to\infty$, and
hence also $N\to\infty$; for $N\ge e$ we have
$L\le2\ln N$. Therefore, the conditions of the corollary imply \eqref{eq:anyp} for
all sufficiently large $n$ and $128L^2/p^2\le512\ln^2N/p^2=o(\alpha)$; it remains to
apply Corollary~\ref{sled:quant}. \hfill$\square$

\emph{Corollary~\ref{sled:lambda}.} Set
$\varepsilon_n:=\bigl(\ln N/(p(\beta-\gamma))\bigr)^{1/2}$; by assumption,
$\varepsilon_n\to0$, and $\ln N=\varepsilon_n\cdot\varepsilon_np(\beta-\gamma)
=o\bigl(\varepsilon_np(\beta-\gamma)\bigr)$. Set
$k:=\lceil(1+\varepsilon_n)\lambda\alpha\rceil$; then $k-\alpha\ge\varepsilon_n\lambda\alpha$,
and for large $n$ we have $k\le2\lambda\alpha$ and $k\to\infty$ (since
$\lambda\alpha\to\infty$). It suffices to show that $\PP(\alpha(G_p)\ge k)\to0$
(for $k>N$ this probability is zero). By Theorem~\ref{th:core} and
the inequalities $(1-p)^x\le e^{-px}$, $\bC Nk\le e^{kL}$, it does not exceed
\begin{displaymath}
\exp\Bigl\{k\Bigl(L+\frac p2-\frac{pk}{2\beta}\Bigr)\Bigr\}+
\exp\bigl\{kL-p(k-\alpha)(\beta-\gamma)\bigr\}.
\end{displaymath}
Since $N\to\infty$ and $p\lambda\alpha/\beta\gg\ln N$, for large $n$
we have $pk/(2\beta)\ge p\lambda\alpha/(2\beta)\ge2L$, and the first exponent does not
exceed $k(L+1/2-2L)\le-k/2$. Further, $p(k-\alpha)(\beta-\gamma)\ge
\varepsilon_n\lambda\alpha\,p(\beta-\gamma)\ge\frac k2\,\varepsilon_np(\beta-\gamma)$,
and $\varepsilon_np(\beta-\gamma)\gg\ln N$, so that for large $n$ we have
$\varepsilon_np(\beta-\gamma)\ge3L$ and the second exponent does not exceed
$kL-\frac32kL=-kL/2$. Both terms tend to zero, since $k\to\infty$.
Thus, a.a.s.\ $\alpha(G_p)<k\le(1+\varepsilon_n)\lambda\alpha+1=(1+o(1))\lambda\alpha$.
\hfill$\square$

\emph{Corollary~\ref{sled:uncond}.} Set $q:=\ln\frac1{1-p}$ and
$\nu:=N$ in the first case, $\nu:=N/\alpha$ in the second; in both cases
$\nu\to\infty$ and $q=o(\ln\nu)$. Fix $\varepsilon\in(0,1)$ and set
$k:=\lceil(1+\varepsilon)\cdot2\alpha\ln\nu/q\rceil$. Since
$\ln\nu/q\to\infty$, for large $n$ we have
$k\le(1+2\varepsilon)\cdot2\alpha\ln\nu/q$ and $k\ge2\alpha\ln\nu/q\ge\alpha$, and it suffices
to prove that $\PP(\alpha(G_p)\ge k)\to0$ (for $k>N$ this probability is
zero). By Theorem~\ref{th:core} with $\beta=\alpha$ and the inequality
$\bC Nk\le(eN/k)^k$
\begin{displaymath}
\PP(\alpha(G_p)\ge k)\ \le\ \exp\Bigl\{k\Bigl(1+\ln\frac Nk+\frac q2
-\frac{qk}{2\alpha}\Bigr)\Bigr\}.
\end{displaymath}
Here $\ln(N/k)\le\ln\nu$ (in the first case because $k\ge1$, in the
second because $k\ge\alpha$) and $qk/(2\alpha)\ge(1+\varepsilon)\ln\nu$;
therefore, the exponent does not exceed $k(1+q/2-\varepsilon\ln\nu)\to-\infty$,
since $q=o(\ln\nu)$ and $\ln\nu\to\infty$. Thus, a.a.s.\
$\alpha(G_p)<k\le(1+2\varepsilon)\cdot2\alpha\log_{1/(1-p)}\nu$, and $\varepsilon$
is arbitrary. \hfill$\square$

\subsection{Proof of Theorem~\ref*{th:order} and Corollary~\ref*{sled:order}}
\label{pf:order}

\emph{(a)} Number the vertices: $V=\{v_1,\dots,v_N\}$. We construct sets
$J_0=\varnothing\subseteq J_1\subseteq\dots\subseteq J_N$: at step $t$
we set $J_t:=J_{t-1}\cup\{v_t\}$ if $|J_{t-1}|<k$ and the vertex $v_t$ is
adjacent in $G_p$ to no vertex of $J_{t-1}$, and $J_t:=J_{t-1}$ otherwise.
The set $J_N$ is independent in $G_p$. At step $t$, only edges
of the graph $G$ between $v_t$ and $J_{t-1}$ are used, and these were not
examined at the previous steps; therefore, for any outcome of the previous steps, if
$|J_{t-1}|<k$, the conditional probability of adding the vertex $v_t$ equals $(1-p)^{h}$, where
$h\le|J_{t-1}|\le k-1$ is the number of neighbors of $v_t$ in $J_{t-1}$, and is therefore at
least $q:=(1-p)^{k-1}$. Realizing the choice at step $t$ by means of an independent
random variable $U_t$ uniformly distributed on $[0,1]$ so that
$v_t$ is certainly added when $U_t\le q$ (if $|J_{t-1}|<k$), we obtain
that $|J_N|\ge\min\{k,X\}$, where $X=\sum_t\mathbf 1\{U_t\le q\}$ has
the binomial distribution with parameters $N$, $q$. By assumption, $\EE X=Nq\ge2k$;
therefore, by Chernoff's inequality (see \cite[Theorem~2.1]{JLR})
\begin{displaymath}
\PP(\alpha(G_p)<k)\le\PP(X<k)\le\PP\Bigl(X\le\frac{\EE X}2\Bigr)\le
e^{-\EE X/8}=\exp\Bigl\{-\frac{N(1-p)^{k-1}}8\Bigr\}\le e^{-k/4}.
\end{displaymath}

\emph{(b)} We may assume that $M\le N$. Since $M\ge16\alpha$ and $\gamma<M/2$,
by Theorem~\ref{th:turan} every $M$-element set contains at least
$M^2/16$ edges of the graph $G$. Therefore,
\begin{displaymath}
\PP(\alpha(G_p)\ge M)\ \le\ \bC NM e^{-pM^2/16}\ \le\
\exp\Bigl\{M\Bigl(1+\ln\frac NM-\frac{pM}{16}\Bigr)\Bigr\}\ \le\ e^{-M}.
\end{displaymath}

\emph{Corollary~\ref{sled:order}.} Let $\gamma\le K\alpha$ for large $n$,
$C:=2K+32$, and $M:=\lceil Cf\alpha\rceil$. Then $M\ge32\alpha\ge16$, $M\ge
(2K+32)\alpha\ge2\gamma+1$, and, since $N/M\le N/\alpha$ and $\ln(N/\alpha)\ge2$
for large $n$,
\begin{displaymath}
pM\ \ge\ Cfp\alpha\ =\ C\ln\frac N\alpha\ \ge\ 32\ln\frac N\alpha\ \ge\
16\Bigl(2+\ln\frac NM\Bigr).
\end{displaymath}
By (b), we have $\PP(\alpha(G_p)\ge M)\le e^{-M}\to0$, since $M\ge f\alpha\ge
\ln(N/\alpha)\to\infty$.

For the lower bound, fix $\eta\in(0,1)$ and set
$k:=\lceil(1-\eta)f\alpha\rceil$; then $k\to\infty$ (since
$f\alpha=\ln(N/\alpha)/p\to\infty$) and $p(k-1)\le(1-\eta)\ln(N/\alpha)$.
For large $n$ we have $p\le\eta/2$ and $\eta f\alpha\ge1$, i.e., $k\le f\alpha$;
since $1-p\ge e^{-p-p^2}$ for $p\le1/2$,
\begin{displaymath}
N(1-p)^{k-1}\ \ge\ Ne^{-(1+p)p(k-1)}\ \ge\ N\Bigl(\frac\alpha N\Bigr)^{(1+p)(1-\eta)}
\ \ge\ \alpha\Bigl(\frac N\alpha\Bigr)^{\eta/2}\ \ge\ 2f\alpha\ \ge\ 2k
\end{displaymath}
for large $n$: here we used the inequalities $(1+p)(1-\eta)\le1-\eta/2$ and
$\frac\eta2\ln\frac N\alpha\ge\ln(2f)$, the second of which holds for large $n$,
since $\ln f=o(\ln(N/\alpha))$. By (a), we have $\PP(\alpha(G_p)<k)\le e^{-k/4}\to0$,
and $\alpha(G_p)\ge k\ge(1-\eta)f\alpha$. The condition
$\gamma=O(\alpha)$ was not used in this part. \hfill$\square$

\subsection{Proof of Theorem~\ref*{th:asthr}}
\label{pf:asthr}

\emph{(ii)} By Remark~\ref{zam:monot}, it suffices to consider $p=cp^*$. Then
$p\to0$ (since $\ln(N/\alpha)=o(\alpha)$) and $f=\ln(N/\alpha)/(p\alpha)=1/c$,
and the statement follows from the lower bound in Corollary~\ref{sled:order} (which
does not use \eqref{eq:struct}).

\emph{(i)} By Remark~\ref{zam:monot}, it suffices to consider
$p=192\lambda/(\varepsilon\alpha)$, where $\lambda:=\ln(eN/\alpha)\ge1$; we may
assume that $p\le1$. Set $M:=\lfloor(1+\varepsilon)\alpha\rfloor+1$; the event
$\alpha(G_p)>(1+\varepsilon)\alpha$ means that $G_p$ contains an independent
set $D$ of cardinality $M$. Let $d:=\stdeg D$; clearly, $d\le\alpha$.

\emph{Case A: $d\le M/16$.} Since $\gamma<\alpha/2<M/2$,
$\psi\le\alpha/16<M/16$, and $M\ge16$ for large $n$, by
Theorem~\ref{th:turan} we have $e(D)\ge M^2/16$. Further,
$pM/16\ge p\alpha/16=12\lambda/\varepsilon\ge12\lambda$; therefore, the sum of the
probabilities over all such $D$ does not exceed
\begin{displaymath}
\bC NMe^{-pM^2/16}\le\exp\Bigl\{M\Bigl(1+\ln\frac NM-12\lambda\Bigr)\Bigr\}
\le e^{-11M}\to0 ,
\end{displaymath}
since $1+\ln(N/M)\le\lambda$.

\emph{Case B: $d>M/16$.} Choose $I\in\I$ such that $|D\cap I|=d$, and set
$j:=|D\setminus I|=M-d\ge M-\alpha>\varepsilon\alpha$. Each vertex
$u\in D\setminus I$ is adjacent to at least $\mu$ vertices of $I$, of which
at most $\alpha-d$ lie outside $D$; hence $u$ is adjacent to at least
$\mu-\alpha+d=d-\delta\alpha\ge d/2$ vertices of $D\cap I$ (since
$\delta\alpha\le\alpha/32<d/2$), and $d/2>M/32>\alpha/32$. The edges found for
distinct $u$ are distinct; therefore, $e(D)\ge j\alpha/32$. The set $D$
is determined by the set $I$ ($m$ possibilities), the set $D\cap I$ (at most
$2^\alpha$ possibilities), and the set $D\setminus I$ (at most $\bC{N-\alpha}j\le
(eN/(\varepsilon\alpha))^j$ possibilities). Therefore, the sum of the probabilities over all
such $D$ does not exceed
\begin{displaymath}
(\alpha+1)\,m\,2^\alpha\max_{j>\varepsilon\alpha}\exp\Bigl\{j\Bigl(1+\ln
\frac N{\varepsilon\alpha}-\frac{p\alpha}{32}\Bigr)\Bigr\}.
\end{displaymath}
By \eqref{eq:lnm}, we have $\ln m\le\alpha\lambda$, and $\ln\frac1\varepsilon\le
\frac1\varepsilon$; since $\lambda\ge1$ and $\varepsilon\le1$,
\begin{displaymath}
1+\ln\frac N{\varepsilon\alpha}+\frac{\ln m+\alpha(1+\ln2)}{\varepsilon\alpha}
\ \le\ \lambda+\frac1\varepsilon+\frac{\lambda+1+\ln2}\varepsilon\ \le\
\frac{5\lambda}\varepsilon\ \le\ \frac{6\lambda}{\varepsilon}=\frac{p\alpha}{32}.
\end{displaymath}
Consequently, for $j>\varepsilon\alpha$ the exponent does not exceed
$-j(\ln m+\alpha+\alpha\ln2)/(\varepsilon\alpha)\le-\ln m-\alpha-\alpha\ln2$, and
the whole sum does not exceed $(\alpha+1)e^{-\alpha}\to0$.

Thus, $\PP(\alpha(G_p)>(1+\varepsilon)\alpha)\to0$. If $p/p^*\to\infty$,
then, since $\lambda=(1+o(1))\ln(N/\alpha)$, for every fixed
$\varepsilon$ and large $n$ we have $p\ge192\lambda/(\varepsilon\alpha)$,
whence $\alpha(G_p)=(1+o(1))\alpha$ a.a.s. \hfill$\square$

\subsection{Proof of Theorem~\ref*{th:exthr}(i)}
\label{pf:exthr1}

The general scheme of the proof is the same as in \cite{BNR,OgRai,PR2}. By
Remark~\ref{zam:monot}, it suffices to consider $p=p_1=64\Lambda/\alpha$
(for $p_1>1$ the statement is trivial). For large $n$, conditions \eqref{eq:struct2} hold, $\alpha\ge e^7$, and
$\Lambda\ge56$. Note that $\ln(N-\alpha)\le\Lambda$ and $\ln m\le\Lambda$.

For $i\ge1$, let $X_i$ be the number of sets $D$ that are independent in $G_p$ and have
cardinality $\alpha$ and degree $\stdeg D=\alpha-i$, and let $Y$ be the number of pairs
$(I,v)$, $I\in\I$, $v\in V\setminus I$, for which the set $I\cup\{v\}$ is
independent in $G_p$. If $Y=0$ and $X_i=0$ for all $i\ge1$, then, first,
every set that is independent in $G_p$ and has cardinality $\alpha$ has degree
$\alpha$, i.e., contains some $I\in\I$, and hence coincides with it;
second, there are no sets that are independent in $G_p$ and have cardinality $\alpha+1$: such a
set would contain an independent subset of cardinality $\alpha$, i.e.,
some $I\in\I$, which contradicts the equality $Y=0$. Hence, in this case $\alpha(G_p)=\alpha$ and
$\I(G_p)=\I(G)$. Set $\varepsilon_0:=1/10$ and split the values of $i$ into
three ranges: (1) $1\le i\le\varepsilon_0\alpha$; (2)
$\varepsilon_0\alpha<i\le\frac{15}{16}\alpha$; (3) $i>\frac{15}{16}\alpha$.

\emph{The quantity $Y$.} The number of pairs $(I,v)$ is exactly $e^\Lambda$, and for each of them all
edges between $v$ and $I$ (there are at least $\mu\ge(1-\delta)\alpha$ of them) must be absent from
$G_p$. Therefore,
\begin{displaymath}
\EE Y\le e^\Lambda(1-p)^{(1-\delta)\alpha}\le\exp\bigl\{\Lambda-64(1-\delta)
\Lambda\bigr\}\le e^{-62\Lambda}\to0 .
\end{displaymath}

\emph{Range (1).} For $j\ge i\ge1$, denote by $X_{i,j}$ the number of
inclusion-maximal sets that are independent in $G_p$ and have cardinality
$\alpha-i+j$ and degree $\alpha-i$. If $D$ is counted in $X_i$ with
$i\le\varepsilon_0\alpha$, extend $D$ to an inclusion-maximal
set $\overline D$ independent in $G_p$; then $\stdeg\overline
D=\alpha-i'$, $0\le i'\le i$. For $i'=0$, the set $\overline D$ contains
some $I\in\I$ and is different from it (otherwise $D\subseteq I$, $|D|=\alpha$, and
$D=I$, which contradicts the inequality $i\ge1$), so that $Y>0$. For $i'\ge1$, the set $\overline D$
is counted in $X_{i',j'}$, where $j'=|\overline D|-(\alpha-i')\ge i'$.
Therefore, it suffices to prove that
$\Sigma_1:=\sum_{i=1}^{\varepsilon_0\alpha}\sum_{j\ge i}\EE X_{i,j}\to0$.

Let $D$ be counted in $X_{i,j}$, let $I\in\I$ be a set on which
the degree is attained, $A_1:=I\setminus D$ ($|A_1|=i$), and $A_2:=D\setminus I$
($|A_2|=j$); the set $D$ is uniquely recovered from the triple $(I,A_1,A_2)$.
Each vertex $u\in A_2$ is adjacent to at least $\mu-i$ vertices
of $D\cap I$, and all these edges (distinct for distinct $u$) are absent from
$G_p$. By the maximality of $D$, each vertex $w\in A_1$ is adjacent in $G_p$
to some vertex of $D$; since $D\setminus A_2\subseteq I$ and $I$
is independent, this vertex lies in $A_2$, and the probability of this event is at most $jp$.
For distinct $w\in A_1$, these events are determined by disjoint sets of
edges and are independent of the edges considered above. Hence
\begin{displaymath}
\EE X_{i,j}\le B_{i,j}:=m\,\bC\alpha i\bC{N-\alpha}j(1-p)^{j(\mu-i)}(jp)^i,
\end{displaymath}
where $\mu-i\ge(1-\delta-\varepsilon_0)\alpha\ge0.8\,\alpha$. If
$B_{i,j+1}>0$, then
\begin{displaymath}
\frac{B_{i,j+1}}{B_{i,j}}=\frac{N-\alpha-j}{j+1}(1-p)^{\mu-i}
\Bigl(1+\frac1j\Bigr)^i\ \le\ e\,(N-\alpha)\,e^{-0.8p\alpha}=:\rho,
\end{displaymath}
and $\ln\rho\le1+\Lambda-51.2\Lambda\le-50\Lambda$. Therefore, $\rho\le1/2$ and
$\sum_{j\ge i}B_{i,j}\le2B_{i,i}$. Further,
\begin{displaymath}
B_{i,i}\le m\Bigl[\frac{e\alpha}i\cdot\frac{e(N-\alpha)}i\cdot ip\cdot
e^{-0.8p\alpha}\Bigr]^i\le mT^i,\qquad T:=e^2(N-\alpha)\,p\alpha\,
e^{-0.8p\alpha},
\end{displaymath}
and $\ln T\le2+\Lambda+\ln(64\Lambda)-51.2\Lambda\le-50\Lambda$ (the last
inequality is equivalent to $2+\ln(64\Lambda)\le0.2\Lambda$ and holds for
$\Lambda\ge56$). Consequently, $T\le1/2$ and
\begin{displaymath}
\Sigma_1\le2m\sum_{i\ge1}T^i\le4mT\le4e^{\Lambda-50\Lambda}\to0 .
\end{displaymath}

\emph{Range (2).} Let $D$ be counted in $X_i$, and let $I\in\I$ be a
set on which the degree $d:=\alpha-i\ge\alpha/16$ is attained. Then
$|I\setminus D|=|D\setminus I|=i$, and each vertex of $D\setminus I$ is adjacent to
at least $\mu-i=d-\delta\alpha\ge(1-16\delta)d\ge\frac34d$ vertices of
$D\cap I$. Since $\frac34pd\ge\frac{3p\alpha}{64}=3\Lambda$ and
$i>\varepsilon_0\alpha$,
\begin{displaymath}
\EE X_i\le m\bC\alpha i\bC{N-\alpha}i(1-p)^{\frac34id}\le
m\Bigl[\frac{e\alpha}i\cdot\frac{e(N-\alpha)}ie^{-3\Lambda}\Bigr]^i\le
m\Bigl[\frac{e^2}{\varepsilon_0^2}\cdot\frac{N-\alpha}{\alpha}\,
e^{-3\Lambda}\Bigr]^i .
\end{displaymath}
The logarithm of the expression in square brackets does not exceed
$2+2\ln10+\Lambda-\ln\alpha-3\Lambda\le-2\Lambda$ (since $\ln\alpha\ge7>
2+2\ln10$). Therefore,
\begin{displaymath}
\sum_{\varepsilon_0\alpha<i\le\frac{15}{16}\alpha}\PP(X_i>0)\le
\alpha\,e^{\Lambda}e^{-2\varepsilon_0\alpha\Lambda}
=\exp\{\ln\alpha+\Lambda-0.2\,\alpha\Lambda\}\to0 .
\end{displaymath}

\emph{Range (3).} Here $\stdeg D<\alpha/16$, and by Theorem~\ref{th:turan}
(with $M=\alpha$) we have $e(D)\ge\alpha^2/16$. Since $N/\alpha\le1+(N-\alpha)\le
2(N-\alpha)$, we have $\ln(N/\alpha)\le\Lambda+\ln2$ and
\begin{displaymath}
\sum_{i>\frac{15}{16}\alpha}\PP(X_i>0)\le\bC N\alpha e^{-p\alpha^2/16}\le
\exp\bigl\{\alpha\bigl(1+\ln2+\Lambda-4\Lambda\bigr)\bigr\}\to0 .
\end{displaymath}
Thus, $\PP\bigl(Y>0\text{ or }\exists\,i\ge1\colon X_i>0\bigr)\to0$, and
statement (i) is proved. \hfill$\square$

\subsection{Proof of Theorem~\ref*{th:exthr}(ii)}
\label{pf:exthr2}

By Remark~\ref{zam:monot}, it suffices to consider $p=p_2=(1-\varepsilon)
\Lambda/\alpha$; since $\Lambda^2=o(\alpha)$, we have $p\to0$, in particular
$p\le1/2$ for large $n$. For a pair $\pi=(I,v)$, $I\in\I$, $v\in V\setminus I$,
let $E_\pi$ be the set of edges of the graph $G$ between $v$ and $I$, $g_\pi:=|E_\pi|\in
[\mu,\alpha]$, and let $\xi_\pi$ be the indicator of the event that all edges of $E_\pi$
are absent from $G_p$, i.e., that the set $I\cup\{v\}$ is independent in $G_p$.
Set $Y:=\sum_\pi\xi_\pi$; if $Y>0$, then $\alpha(G_p)\ge\alpha+1$.
The number of pairs equals $K:=e^\Lambda$, and since $1-p\ge e^{-p-p^2}$,
\begin{displaymath}
\EE Y\ \ge\ K(1-p)^\alpha\ \ge\ \exp\{\Lambda-p\alpha-p^2\alpha\}
=\exp\Bigl\{\varepsilon\Lambda-\frac{(1-\varepsilon)^2\Lambda^2}{\alpha}\Bigr\}
\to\infty.
\end{displaymath}
Further, $\EE\bigl(Y(Y-1)\bigr)=\sum_{\pi\ne\pi'}(1-p)^{g_\pi+g_{\pi'}-|E_\pi\cap
E_{\pi'}|}$. Let $\pi=(I,v)\ne\pi'=(I',v')$.

1) $v\ne v'$. The sets $E_\pi$ and $E_{\pi'}$ can have in common only the edge
$\{v,v'\}$; therefore, $(1-p)^{-|E_\pi\cap E_{\pi'}|}\le(1-p)^{-1}\le1+2p$, and the contribution of
such pairs does not exceed $(1+2p)(\EE Y)^2$.

2) $v=v'$, $I\ne I'$. The common edges have the form $\{v,w\}$, $w\in I\cap I'$; there are at
most $\theta$ of them, so that the term is at most $(1-p)^{2\mu-\theta}$. The number of such
ordered pairs is at most $m^2(N-\alpha)=K^2/(N-\alpha)$ (we choose $I$, $I'$, and
then $v\notin I$), and $(\EE Y)^2\ge K^2(1-p)^{2\alpha}$. Therefore, their contribution does not
exceed
\begin{displaymath}
(\EE Y)^2\,\frac{(1-p)^{-(2\delta\alpha+\theta)}}{N-\alpha}\le(\EE Y)^2
\exp\bigl\{2p(2\delta\alpha+\theta)-\ln(N-\alpha)\bigr\}\le(\EE Y)^2
\exp\Bigl\{4\Lambda\Bigl(\delta+\frac\theta\alpha\Bigr)-\ln(N-\alpha)\Bigr\},
\end{displaymath}
and this is $o\bigl((\EE Y)^2\bigr)$ by \eqref{eq:second}; here we used
the inequalities $(1-p)^{-1}\le e^{2p}$ and $p\le\Lambda/\alpha$.

Thus, $\EE\bigl(Y(Y-1)\bigr)\le(1+o(1))(\EE Y)^2$, whence
$\DD Y\le o\bigl((\EE Y)^2\bigr)+\EE Y$, and by Chebyshev's inequality
$\PP(Y=0)\le\DD Y/(\EE Y)^2\to0$. \hfill$\square$

\section{Concluding remarks}
\label{sec:concl}

1. The constants in Corollary~\ref{sled:quant} and Theorems~\ref{th:turan}--\ref{th:exthr}
are not optimal. In Theorems~\ref{th:asthr} and~\ref{th:exthr} they are due
mainly to the constants $\frac34$ and $\frac1{16}$ in Lemma~\ref{lem:bad}, and in
Corollary~\ref{sled:quant} to the choice of the parameter $\beta$ in the proof. As noted in \cite[Section~6]{PR2} using the graphs $G(n,r,<s)$ as an example, the sharp threshold (up to a factor $1+o(1)$)
apparently cannot be obtained by this method: this requires more delicate considerations,
which so far have been developed only for particular families \cite{BNR,BKL}.

2. The conditions of Theorems~\ref{th:asthr} and~\ref{th:exthr} are expressed (apart from $N$ and
$\alpha$) in terms of five quantities $\gamma,\psi,\mu,\theta,m$, and all of them, as the
examples show, are estimated in an elementary way if the structure of the
independent sets of maximum cardinality is known. As a rule, the hardest part is
estimating the quantity $\psi$: for the graphs $G(n,r,<s)$, in one of the two cases,
\cite{PR2} invokes the Ahlswede--Khachatrian theorem \cite{AK2} for this purpose.

3. The graphs $G(n,r,s)$ (vertices are adjacent if they intersect in exactly $s$ elements)
with $s\ge1$ and constant $r$ do not satisfy the conditions of our theorems: for two adjacent vertices $x,y$, the common non-neighbors
include, in particular, all sets disjoint from $x\cup y$, and there are
asymptotically as many of them as there are vertices in total (for $s=1$ this is noted in
\cite[Subsection~3.3.10]{RaiSin}), i.e., $\gamma\asymp N$. Accordingly, the threshold for asymptotic stability
$\ln n/n^{r-2s-1}$ found in \cite{Pyad19} for constant $r>2s+1$
has a form different from $p^*$. The threshold for exact stability for such graphs with $r>2s+2$ is studied in
\cite{Kosh}. However, the unconditional bounds of
Corollary~\ref{sled:uncond} and Remark~\ref{zam:lower} apply to the graphs $G(n,r,s)$. For example, for
$G(n,3,1)$ (here $\alpha\sim n$, see \cite[Subsection~4.2]{Rai26}), for
$\ln\frac1p=o(\ln n)$ and $\ln\frac1{1-p}=o(\ln n)$, Corollary~\ref{sled:uncond}
gives a.a.s.\ $\alpha(G_p)\le(1+o(1))\cdot4n\log_{1/(1-p)}n$, and
Remark~\ref{zam:lower}, applied to $f\sim n$ pairwise nonadjacent
complete subgraphs on $w\sim n/\ln n$ vertices from \cite[Subsection~4.2]{Rai26},
gives a.a.s.\ $\alpha(G_p)\ge(1-o(1))\,n\log_{1/(1-p)}n$; for $p=1/2$, the exact
asymptotic value $2n\log_2n$ was found in \cite{Pyad16} (for the qualitative picture for
$r\le2s+1$, see \cite[Theorem~10]{RaiSin}). We do not know general conditions giving the threshold
for asymptotic stability for the graphs $G(n,r,s)$ with $r>2s+1$.

4. In the second example from \cite{Rai17}, constructed from a Sidon 4-set in an
abelian group, the quantity $\gamma$ is bounded. To apply
Theorems~\ref{th:asthr} and~\ref{th:exthr} to it, one needs to describe the independent sets
of maximum cardinality of this graph and to bound the quantity $\psi$.


\end{document}